\documentclass[11pt,a4paper]{article}

\usepackage{titlesec}
\usepackage{fancyhdr}
\usepackage{a4wide}
\usepackage{graphicx}
\usepackage{float}
\usepackage{amssymb}
\usepackage{amsmath}
\usepackage{amsthm}
\usepackage{color}
\usepackage{xcolor,colortbl}
\usepackage{mathrsfs}
\usepackage{array}
\usepackage{eucal}
\usepackage{tikz}
\usepackage[T1]{fontenc}
\usepackage{inputenc}
\usepackage[english]{babel}
\usepackage{lmodern}
\usepackage{hyperref}
\usepackage{geometry}
\usepackage{changepage}
\changepage{0pt}{}{}{}{}{0pt}{}{0pt}{10pt}
\usepackage[numbers]{natbib}
\usepackage{hyperref}
\hypersetup{
pdfpagemode=none,
pdftoolbar=true,        
pdfmenubar=true,        
pdffitwindow=false,     
pdfstartview={Fit},    
pdftitle={Trapped modes and reflectionless modes as eigenfunctions of the same  spectral problem},    
pdfauthor={A.-S. Bonnet-Ben Dhia, L. Chesnel, V. Pagneux},     
pdfsubject={},  
pdfcreator={A.-S. Bonnet-Ben Dhia, L. Chesnel, V. Pagneux},   
pdfproducer={A.-S. Bonnet-Ben Dhia, L. Chesnel, V. Pagneux}, 
pdfkeywords={}, 
pdfnewwindow=true,      
colorlinks=true,       
linkcolor=magenta,          
citecolor=red,        
filecolor=cyan,      
urlcolor=blue           
}

\newcommand{\dsp}{\displaystyle}

\newcommand{\om}{\omega}
\newcommand{\Om}{\Omega}
\newcommand{\mrm}[1]{\mathrm{#1}}

\newcommand{\Cplx}{\mathbb{C}}
\newcommand{\N}{\mathbb{N}}
\newcommand{\R}{\mathbb{R}}

\newcommand{\mL}{\mrm{L}}
\newcommand{\mH}{\mrm{H}}

\renewcommand{\ker}{\mrm{ker}\,}

\newcommand{\Indicator}{\rho}
\definecolor{RoyalBlue}{cmyk}{1, 0.50, 0, 0}
\definecolor{Gray}{gray}{0.90}

\newtheorem{theorem}{Theorem}[section]

\newtheorem{definition}{Definition}[section]

\newtheorem{proposition}{Proposition}[section]

\begin{document}

~\vspace{0.0cm}
\begin{center}
{\sc \bf\LARGE  Trapped modes and reflectionless modes as\\[6pt]  eigenfunctions of the same  spectral problem}
\end{center}

\begin{center}
\textsc{Anne-Sophie Bonnet-Ben Dhia}$^1$, \textsc{Lucas Chesnel}$^2$, \textsc{Vincent Pagneux}$^3$\\[16pt]
\begin{minipage}{0.96\textwidth}
{\small

$^1$ Laboratoire Poems, CNRS/ENSTA/INRIA, Ensta ParisTech, Universit\'e Paris-Saclay, 828, Boulevard des Mar\'echaux, 91762 Palaiseau, France;\\
$^2$ INRIA/Centre de math\'ematiques appliqu\'ees, \'Ecole Polytechnique, Universit\'e Paris-Saclay, Route de Saclay, 91128 Palaiseau, France;\\
$^3$ Laboratoire d'Acoustique de l'Universit\'e du Maine, Av. Olivier Messiaen, 72085 Le Mans, France.\\[10pt]
E-mails: \texttt{bonnet@ensta.fr},  \texttt{lucas.chesnel@inria.fr}, \texttt{vincent.pagneux@univ-lemans.fr}\\[-12pt]
\begin{center}
-- \today --
\end{center}
}
\end{minipage}
\end{center}
\vspace{0.4cm}

\noindent\textbf{Abstract.} We consider the reflection-transmission problem in a waveguide with obstacle. At certain frequencies, for some incident waves, intensity is perfectly transmitted and the reflected field decays exponentially at infinity. In this work, we show that such reflectionless modes can be characterized as eigenfunctions of an original non-selfadjoint spectral problem. In order to select ingoing waves on one side of the obstacle and outgoing  waves on the other side, we use complex scalings (or Perfectly Matched Layers) with imaginary parts of different signs. We prove that the real eigenvalues of the obtained spectrum correspond either to trapped modes (or bound states in the continuum) or to reflectionless modes. Interestingly, complex eigenvalues also contain useful information on weak reflection cases. When the geometry has certain symmetries, the new spectral problem enters the class of $\mathcal{PT}$-symmetric problems.\\


\section{Introduction}

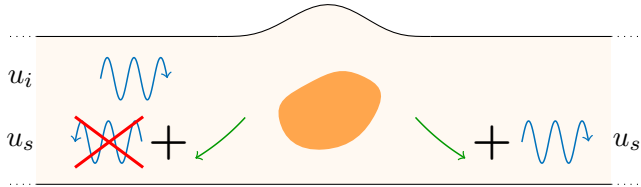
\begin{figure}[!ht]
\centering
\begin{tikzpicture}[scale=1.4]
\draw[fill=orange!5,draw=none](-2.7,0.3) rectangle (2.7,1.7);
\draw[samples=30,domain=-1:1,draw=black,fill=orange!5] plot(\x,{1.7+0.1*(\x+1)^4*(\x-1)^4*(\x+3)});
\begin{scope}[scale=0.7,yshift=0.4cm]
\draw [fill=orange!70,draw=none] plot [smooth cycle, tension=1] coordinates {(-0.6,0.9) (0,0.5) (0.7,1) (0.5,1.5) (-0.2,1.4)};
\end{scope}
\draw (-2.7,0.3)--(2.7,0.3);
\draw (-2.7,1.7)--(-1,1.7);
\draw (1,1.7)--(2.7,1.7);
\draw [dotted](-2.7,0.3)--(-3,0.3);
\draw [dotted](2.7,0.3)--(3,0.3);
\draw [dotted](-2.7,1.7)--(-3,1.7);
\draw [dotted](2.7,1.7)--(3,1.7);
\node at (-2.85,1.3){ $u_i$};
\node at (-2.85,0.7){$u_s$};
\node at (2.85,0.7){$u_s$};
\begin{scope}[xshift=0.9cm,yshift=0.7cm,scale=0.8]
\begin{scope}[xshift=-0.05cm]
\draw[green!60!black,line width=0.2mm,->] plot[very thick,domain=0:pi/4-0.2,samples=100] (\x,{-0.8+0.02*exp(-(\x-4))});
\end{scope}
\node at (0.85,0){\LARGE $\boldsymbol{+}$};
\begin{scope}[xshift=1.2cm]
\draw[RoyalBlue,line width=0.2mm,->] plot[domain=0:pi/4,samples=100] (\x,{0.25*sin(20*\x r)});
\end{scope}
\end{scope}

\begin{scope}[xshift=-2.5cm,yshift=0.7cm,scale=0.8]
\begin{scope}[xshift=0cm]
\draw[RoyalBlue,line width=0.2mm,<-] plot[domain=0:pi/4,samples=100] (0.2+\x,{0.25*sin(20*\x r)});
\node at (1.3,0){\LARGE $\boldsymbol{+}$};
\draw[red,line width=0.4mm] (0.2,-0.3)--(1,0.3);
\draw[red,line width=0.4mm] (1,-0.3)--(0.2,0.3);
\end{scope}
\begin{scope}[xshift=2.5cm,yshift=-0.8cm]
\draw[green!60!black,line width=0.2mm,<-] plot[very thick,domain=8-pi/4+0.2:8,samples=100] (\x-8.3,{0.02*exp((\x-4))});
\end{scope}
\end{scope}
\begin{scope}[xshift=-2.1cm,yshift=1.3cm,scale=0.8]
\draw[RoyalBlue,line width=0.2mm,->] plot[domain=0:pi/4,samples=100] (\x,{0.25*sin(20*\x r)});
\end{scope}
\end{tikzpicture}
\caption{Schematic picture of a reflectionless mode. The propagating wave (blue) is not reflected. The backscattered field is purely evanescent (green). \label{SchematicPicture}
} 
\label{figIntro}
\end{figure}

In bounded domains, eigenvalues and corresponding eigenmodes completely determine the solutions to wave equations. In open systems, it is necessary to take into account also \textit{complex resonances}, sometimes known as \textit{leaky modes} or \textit{quasi normal modes} in literature \cite{aguilar1971class,balslev1971spectral,simon1972quadratic,moiseyev1998quantum,AsPV00,sjostrand1991complex,zworski1999resonances,zworski2011lectures,moiseyev2011non}.  Complex resonances provide the backbone of phenomena of wave propagation. For example, they indicate if there is a rapid variation of the scattering coefficients. We emphasize that real eigenvalues and complex resonances are intrinsic and universal objects that can be computed using different approaches. Among them, let us mention methods based on complex scaling or \textit{Perfectly Matched Layers} (PMLs) \cite{Bera94} or techniques relying on Hardy space infinite elements \cite{HoNa09,HoNa15}. 
However complex resonances fail to provide useful information on recurrent questions such as ``how large is the backscattered energy?'' or ``is there good transmission in the system?''. 
Such questions are of tremendous importance in many topics of intense study in wave physics: extraordinary optical transmission \cite{ebbesen1998extraordinary},  
topological states immune to backscattering  \cite{wang2009observation,lu2014topological},
perfect transmission resonances \cite{zhukovsky2010perfect}, 
transmission eigenchannels through disordered media \cite{kim2012maximal,sebbah2017scattering},
reflectionless metamaterials \cite{rahm2008optical} or
metasurfaces \cite{yu2014flat,asadchy2015broadband}.\\
\newline
In this article, we define a new complex spectrum to quantify these aspects.  In particular, our spectrum allows one to identify wavenumbers for which there is an incident field such that the backscattered field is evanescent (Figure \ref{figIntro}). For construction of waveguides where such phenomenon occurs at a prescribed frequency, we refer the reader to \cite{BoNa13,ChNa16,BoCN17,ChNPSu}. Our spectrum also contains usual trapped modes \cite{Urse51,Evan92,EvLV94,DaPa98,LiMc07,Naza10c} called \textit{Bound States in the Continuum} (BSCs or BICs) in quantum mechanics \cite{SaBR06,Mois09,GPRO10,HZSJS16}. To compute this new spectrum, we use complex scaling/PMLs techniques in an original way. Our approach is based on the following basic observation: if for an incident wave, the backscattered field is evanescent, then the total field is ingoing in the input lead and outgoing in the output lead. Therefore, changing the sign of the imaginary part of the usual complex scaling/PML in the input lead should allow one to compute \textit{reflectionless modes} (see the exact definition after (\ref{defEquationDilatation})). We show that this simple idea, inspired by the results obtained in \cite{HeKS11} on a 1D $\mathcal{PT}$-symmetric problem, indeed works. Interestingly, we will observe numerically that complex eigenvalues in this spectrum seem to contain useful information as well, indicating settings where reflection is weak.

\noindent To make the presentation as simple as possible, we stick to a simple 2D wave problem in a waveguide with Neumann boundary conditions. We assume that this waveguide contains a penetrable obstacle (Figure \ref{figsetting}). Everything presented here can be generalized to other types of obstacles, to other kinds of boundary conditions (Dirichlet, ...), to higher dimension ($\mrm{3D}$) and to more complex geometries (one input lead/several output leads).\\
\newline
The paper is organized as follows. In section \ref{SectionSetting}, the governing equations of the scattering problem are presented and the reflectionless case is defined. Section \ref{SectionClassical} is a reminder of the technique of complex scaling allowing one to compute  classical complex resonances. In section \ref{SectionUnusual}, we explain how to use conjugated complex scalings to identify reflectionless frequencies. We provide examples of spectra in section \ref{SectionNumerics} dedicated to numerics and we give some proofs of the theoretical results in section \ref{SectionProofs}. Eventually, section \ref{SectionConclusion} contains some concluding remarks.

\section{Setting}\label{SectionSetting}

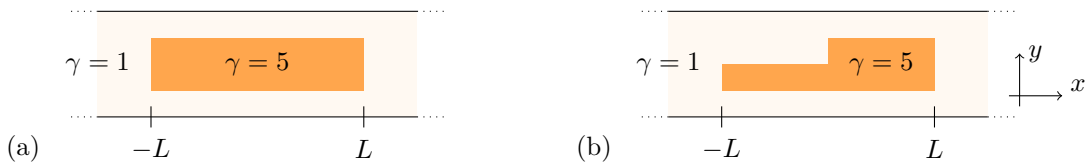
\begin{figure}[!ht]
\centering
\begin{tikzpicture}[scale=1.4]
\draw[fill=orange!5,draw=none](-1.5,0) rectangle (1.5,1);
\draw[fill=orange!70,draw=none](-1,0.25) rectangle (1,0.75);
\draw (-1.5,0)--(1.5,0);
\draw (-1.5,1)--(1.5,1);
\draw [dotted](-1.5,0)--(-1.8,0);
\draw [dotted](1.5,0)--(1.8,0);
\draw [dotted](-1.5,1)--(-1.8,1);
\draw [dotted](1.5,1)--(1.8,1);
\draw (-1,-0.1)--(-1,0.1);
\draw (1,-0.1)--(1,0.1);
\node at (-1,-0.3){\small $-L$};
\node at (1,-0.3){\small $L$};
\node at (0,0.5){\small $\gamma=5$};
\node at (-1.5,0.5){\small $\gamma=1$};
\node at (-2.2,-0.3){\small (a)};
\end{tikzpicture}\qquad\qquad\begin{tikzpicture}[scale=1.4]
\draw[fill=orange!5,draw=none](-1.5,0) rectangle (1.5,1);
\draw[fill=orange!70,draw=none](-1,0.25) rectangle (1,0.75);
\draw[fill=orange!5,draw=none](-1,0.5) rectangle (0,0.75);
\draw (-1.5,0)--(1.5,0);
\draw (-1.5,1)--(1.5,1);
\draw [dotted](-1.5,0)--(-1.8,0);
\draw [dotted](1.5,0)--(1.8,0);
\draw [dotted](-1.5,1)--(-1.8,1);
\draw [dotted](1.5,1)--(1.8,1);
\node at (0.5,0.5){\small $\gamma=5$};
\node at (-1.5,0.5){\small $\gamma=1$};
\begin{scope}[xshift=-1.3cm,yshift=-1cm]
\draw[->] (3,1.2)--(3.5,1.2);
\draw[->] (3.1,1.1)--(3.1,1.6);
\node at (3.65,1.3){\small $x$};
\node at (3.25,1.6){\small $y$};
\end{scope}
\draw (-1,-0.1)--(-1,0.1);
\draw (1,-0.1)--(1,0.1);
\node at (-1,-0.3){\small $-L$};
\node at (1,-0.3){\small $L$};
\node at (-2.2,-0.3){\small (b)};
\end{tikzpicture}\\[-6pt]
\caption{Symmetric (a) and non-symmetric (b) obstacles considered in the numerics below. \label{figsetting}} 
\end{figure}

\noindent We are interested in the propagation of waves in a waveguide with a bounded penetrable obstacle (Figure \ref{figsetting}). We assume that the waveguide coincides with the region $\Om:=\{(x,y)\in\R^2\,|\,0<y<1\}$ and that the propagation of waves is governed by the Helmholtz equation with Neumann boundary conditions 
\begin{equation}\label{PbInitial}
\begin{array}{|rcll}
\Delta u + k^2\gamma u & = & 0 & \mbox{ in }\Om\\[3pt]
 \partial_y u & = & 0  & \mbox{ on }y=0\mbox{ and }y=1.
\end{array}
\end{equation}
Here $u$ can represent for instance the acoustic pressure in a compressible fluid, or the transverse component of the displacement field in an elastic solid, or one component of the electromagnetic  field in a dielectric medium.  $\Delta$ denotes the $\mbox{2D}$ Laplace operator. Moreover $k$ is the wavenumber such that $k=\om/c$ where $c$ denotes  the  waves velocity and $\om$  is the angular frequency, corresponding to a time harmonic dependence in $e^{-i\om t}$. We assume that $\gamma$ is a positive and bounded function such that $\gamma=1$ for $|x|\ge L$ where $L>0$ is given. In other words, the obstacle is located in the region $\Om_L:=\{(x,y)\in\Om\,|\,|x|<L\}$. \\
Pick $k\in(N\pi;(N+1)\pi)$, with $N\in\N:=\{0,1,\dots\}$. For $n\in\N$, set $\beta_n=\sqrt{k^2-n^2\pi^2}$. In this work, the square root is such that if $z=|z|\,e^{i\mrm{arg}(z)}$ with $\mrm{arg}(z)\in[0;2\pi)$, then $\sqrt{z}=\sqrt{|z|}\,e^{i\mrm{arg}(z)/2}$. With this choice, we have $\Im m\,\sqrt{z}\ge0$ for all $z\in\Cplx$. The function 
\begin{equation}\label{defModes}
w_n^{\pm}(x,y)=(2|\beta_n|)^{-1/2}e^{\pm i\beta_n x}\varphi_n(y),
\end{equation}
with $\varphi_0(y)=1$, $\varphi_n(y)=\sqrt{2}\cos(n\pi y)$ for $n\ge1$, solves Problem (\ref{PbInitial}) without obstacle ($\gamma\equiv1$). For $n=0,\dots,N$, the wave $w_n^{\pm}$ propagates along the $(Ox)$ axis from $\mp\infty$ to $\pm\infty$. On the other hand, for $n> N$, $w_n^{\pm}$ is exponentially  growing at $\mp\infty$ and exponentially decaying at $\pm\infty$. For $n=0,\dots,N$, we consider the scattering of the wave $w_n^{+}$ by the obstacle located in $\Om$. It is known that Problem (\ref{PbInitial}) admits a solution $u_n=w_n^{+}+u^s_n$ with the \textit{outgoing} \textit{scattered field} $u^s_n$ written as
\begin{equation}\label{RadiationCondition}
u^s_n=\sum_{p=0}^{+\infty} s^{\pm}_{np}w_p^{\pm}\qquad\mbox{ for }\pm x\ge L
\end{equation}
with $(s^{\pm}_{np})\in\Cplx^{\N}$. The solution $u_n$ is uniquely defined if and only if \textit{Trapped Modes} (TMs) do not exist at the wavenumber $k$. We remind the reader that trapped modes are non zero  functions $u\in \mL^2(\Omega)$ satisfying (\ref{PbInitial}), where $\mL^2(\Omega)$ is the space of square-integrable functions in $\Omega$. We denote by $\mathscr{K}_{\mrm{t}}$ the set of $k^2$ such that TMs exist at the wavenumber $k$. On the other hand, the scattering coefficients $s^{\pm}_{np}$ in (\ref{RadiationCondition}) are always uniquely defined, including for $k^2\in\mathscr{K}_{\mrm{t}}$. In the following, we will be particularly interested in the features of the \textit{reflection matrix}  (whose size, determined by the number of propagative modes, depends on $k$)
\begin{equation}\label{reflection matrix} 
R(k):=(s^{-}_{np})_{0\le n,p\le N}\in\Cplx^{N+1\times N+1}.
\end{equation} 
\begin{definition} Let  $k\in(0;+\infty)\setminus\pi\N$. 
We say that the wavenumber $k$ is reflectionless if $\ker R(k)\ne\{0\}$.
\end{definition} 
\noindent Let us explain this definition. In general, by linearity, for an incident field 
\begin{equation}\label{DefIncidentField}
u_i=\sum_{n=0}^Na_n w_n^{+},\qquad (a_n)_{n=0}^N\in\Cplx^{N+1},
\end{equation}
Problem (\ref{PbInitial}) admits a solution $u$ such that $u=u_i+u_s$ with
\begin{equation}\label{DefScatteredField}
u_s=\sum_{p=0}^{+\infty}b_p^{\pm} w_p^{\pm}\mbox{ for }\pm x\ge L\ \mbox{ and }\  b_p^{\pm}=\sum_{n=0}^Na_{n}s^{\pm}_{np}\in\Cplx.
\end{equation}
The above definition says that, if $k$ is reflectionless, then there are $(a_n)_{n=0}^N\in\Cplx^{N+1}\setminus\{0\}$ such that the $b_p^{-}$ in (\ref{DefScatteredField}) satisfy $b_p^{-}=0$, $p=0,\dots,N$. In other words, the scattered field is exponentially decaying for $x\le -L$. Let us notice finally that the corresponding total field $u=u_i+u_s$ decomposes as
\begin{equation}\label{DefReflectionlessMode}
\begin{array}{ll}
\dsp u=\sum_{n=0}^{N}a_n w_n^{+}+\tilde{u}&\mbox{ for } x\le -L\\[13pt]
\dsp u=\sum_{n=0}^{N}t_n w_n^{+}+\tilde{u}&\mbox{ for } x\ge L\\
\end{array}
\end{equation}
where $t_n=a_n+b_n^+$ and where $\tilde{u}$ decays exponentially for $\pm x\ge L$. In other words, the total field is \textit{ingoing} for $x\le-L$ and \textit{outgoing} for $x\ge L$.\\
\newline
In the following,  we call \textit{Reflectionless Modes} (RMs) the functions $u$ admitting expansion (\ref{DefReflectionlessMode}) and we denote by $\mathscr{K}_{\mrm{r}}$ the set of $k^2$ such that the wavenumber $k$ is reflectionless. The main objective of this article is to explain how to determine directly the set $\mathscr{K}_{\mrm{r}}$ and the corresponding RMs by solving an eigenvalue problem, instead of computing the reflection matrix for all values of $k$.  

\section{Classical complex scaling}\label{SectionClassical}
As a first step, in this section we remind briefly how to use a complex scaling to compute trapped modes. Define the unbounded operator $A$ of $\mL^2(\Om)$ such that
\[
Au=-\cfrac{1}{\gamma}\,\Delta u
\]
with Neumann boundary conditions $\partial_{y}u=0\mbox{ on }y=0\mbox{ and }y=1$. It is known that $A$ is a selfadjoint operator ($\mL^2(\Om)$ is endowed with the inner product $(\gamma\,\cdot,\cdot)_{\mL^2(\Om)}$) whose spectrum $\sigma(A)$ coincides with $[0;+\infty)$. More precisely, we have $\sigma_{\mrm{ess}}(A)=[0;+\infty)$ where $\sigma_{\mrm{ess}}(A)$ denotes the essential spectrum of $A$. By definition, $\sigma_{\mrm{ess}}(A)$ corresponds to the set of $\lambda \in\Cplx$ for which there exists a so-called singular sequence $(u^{(m)})$, that is an orthonormal sequence $(u^{(m)})\in\mL^2(\Omega)^{\N}$ such that $((A-\lambda)u^{(m)})$ converges to 0 in $\mL^2(\Omega)$. 
Besides, $\sigma(A)$ may contain eigenvalues (at most a sequence accumulating at $+\infty$) corresponding to TMs. In order to reveal these eigenvalues which are embedded in $\sigma_{\mrm{ess}}(A)$, one can use a complex change of variables. For $0<\theta<\pi/2$, set $\eta=e^{i\theta}$ and define the function $\mathcal{I}_{\theta}:\R\to\Cplx $ such that 
\begin{equation}\label{defDilatationx}
\mathcal{I}_{\theta}(x)=\begin{array}{|ll}
-L+(x+L)\,\eta&\mbox{ for }x\le- L\\
x&\mbox{ for }|x|< L\\
+ L+(x-L)\,\eta&\mbox{ for } x\ge L.
\end{array}
\end{equation}
For the sake of simplicity, we  will use abusively the same notation $\mathcal{I}_{\theta}$ for the following map:
$\{\Om \to \Cplx\times(0;1), \;
(x,y) \mapsto (\mathcal{I}_{\theta}(x),y)\}$.
Note that with this definition, the left inverse $\mathcal{I}_{\theta}^{-1}$ of $\mathcal{I}_{\theta}$, acting from $\mathcal{I}_{\theta}(\Om)$ to $\Om$, is equal to  $\mathcal{I}_{-\theta}$. One can easily check that for all $n\geq 0$, 
$w_n^+\circ \mathcal{I}_{\theta}$ is exponentially decaying for $x\ge L$, while $w_n^-\circ \mathcal{I}_{\theta}$ is exponentially decaying for $x\le-L$. As a consequence, defining from expansion (\ref{DefScatteredField})  the function $v_{\theta}=u_s\circ \mathcal{I}_{\theta}$, one has $v_{\theta}=u_s$ for $|x|< L$ and $v_{\theta}\in \mL^2(\Omega)$ (which is in general not true for $u_s$). Moreover $v_{\theta}$ satisfies the following equation in $\Om$:
\begin{equation}\label{defEquationDilatation}
\alpha_{\theta}\frac{\partial}{\partial x}\Big(\alpha_{\theta}\frac{\partial v_{\theta}}{\partial x}\Big)+\frac{\partial^2 v_{\theta}}{\partial y^2}+k^2\gamma v_{\theta}=k^2(1-\gamma)u_i
\end{equation}
with $\alpha_{\theta}(x)=1$ for $|x|<L$ and $\alpha_{\theta}(x)=\eta^{-1}=\overline{\eta}$ for $\pm x\ge L$. In particular, for a TM, $v_{\theta}$ solves (\ref{defEquationDilatation}) with $u_i=0$. This leads us to consider the unbounded operator $A_{\theta}$ of $\mL^2(\Om)$ such that
\begin{equation}\label{defOpClassicalPMLs}
A_{\theta}v_{\theta}=-\cfrac{1}{\gamma}\left(\alpha_{\theta}\frac{\partial}{\partial x}\Big(\alpha_{\theta}\frac{\partial v_{\theta}}{\partial x}\Big)+\frac{\partial^2 v_{\theta}}{\partial y^2}\right)
\end{equation}
again with homogeneous Neumann boundary conditions. Since $\alpha_{\theta}$ is complex valued, the operator $A_{\theta}$ is no longer selfadjoint. However, we use the same definition as above for $\sigma_{\mrm{ess}}(A_{\theta})$, which is licit for this operator. We recall below the main spectral properties of $A_{\theta}$ \cite{Simo78}:
\begin{theorem}\label{thmUsualPMLs}
i) There holds
\begin{equation}\label{essClassicalPMLs}
\sigma_{\mrm{ess}}(A_{\theta})=\bigcup_{n\in\N,\,t\ge0}\{n^2\pi^2+te^{-2i\theta}\}.
\end{equation}
ii)  The spectrum of $A_{\theta}$ satisfies $\sigma(A_{\theta})\subset\mathscr{R}_{\mrm{\theta}}^-$ with 
\[
\mathscr{R}_{\mrm{\theta}}^-:=\{z\in\Cplx\,|\,-2\theta\le\mrm{arg}(z)\le 0\}.
\]
iii)  $\sigma (A_{\theta})\setminus \sigma_{\mrm{ess}}(A_{\theta})$ is discrete and contains only eigenvalues of finite multiplicity.\\[3pt]
iv) Assume that $k^2\in\sigma(A_{\theta})\setminus\sigma_{\mrm{ess}}(A_{\theta})$. Then $k^2$ is real if and only if $k^2\in\mathscr{K}_{\mrm{t}}$. Moreover if $v_{\theta} $ is an eigenfunction associated to $k^2$ such that $\Im m\,k^2<0$, then $v_{\theta}\circ\mathcal{I}_{-\theta}$ is a solution of the original problem (\ref{PbInitial}) whose amplitude is exponentially growing at $+\infty$ or at $-\infty$.
\end{theorem} 
\noindent The interesting point is that now TMs correspond to isolated eigenvalues of $A_{\theta}$, and as such, they can be computed numerically as illustrated below. Note that the elements $k^2$ of $\sigma(A_{\theta})\setminus\sigma_{\mrm{ess}}(A_{\theta})$ such that $\Im m\,k^2<0$, if they exist, correspond to complex resonances (quasi normal modes). Let us point out that the complex scaling is just a technique to reveal them. Indeed complex resonances are intrinsic objects defined as the poles of the meromorphic extension from $\{z\in\Cplx\,|\,\Im m\,z>0\}$ to $\{z\in\Cplx\,|\,\Im m\,z\le0\}$ of the operator valued map $z\mapsto (\Delta +z\gamma)^{-1}$. For more details, we refer the reader to \cite{AsPV00}.

\section{Conjugated complex scaling}\label{SectionUnusual}

Now, we show that replacing the classical complex scaling by an unusual \textit{conjugated} complex scaling, and proceeding  as in the previous section, we can define a new complex spectrum which contains the reflectionless values  $k^2\in \mathscr{K}_{\mrm{r}}$ we are interested in. We define the map $\mathcal{J}_{\theta}:\Om\to\Cplx\times(0;1)$ using the following  complex change of variables  
\begin{equation}\label{defDilatationConjugated}
\mathcal{J}_{\theta}(x)=\begin{array}{|ll}
-L+(x+L)\,\overline{\eta}&\mbox{ for }x\le- L\\
x&\mbox{ for }|x|< L\\
+ L+(x-L)\,\eta&\mbox{ for } x\ge L,
\end{array}
\end{equation}
with again $\eta=e^{i\theta}$ ($0<\theta<\pi/2$). Note the important difference in the definitions of $\mathcal{J}_{\theta}$ and $\mathcal{I}_{\theta}$ for $x\le-L$: $\eta$ has been replaced by the conjugated parameter $\overline{\eta}$ to select the ingoing modes instead of the outgoing ones in accordance with (\ref{DefReflectionlessMode}). Now, if $u$ is a RM associated to  $k^2\in \mathscr{K}_{\mrm{r}}$, setting  $w_{\theta}=u\circ \mathcal{J}_{\theta}$, one has $w_{\theta}=u$ for $|x|< L$ and $w_{\theta}\in \mL^2(\Omega)$ (which is not the case for $u$). The function $w_{\theta}$ satisfies the following equation in $\Om$:
\begin{equation}\label{DefPbFortCPMLs}
\beta_{\theta}\frac{\partial}{\partial x}\Big(\beta_{\theta}\frac{\partial w_{\theta}}{\partial x}\Big)+\frac{\partial^2 w_{\theta}}{\partial y^2}+k^2\gamma w_{\theta}=0
\end{equation}
with $\beta_{\theta}(x)=1$ for $|x|< L$, $\beta_{\theta}(x)=\eta$ for $x\le-L$ and $\beta_{\theta}(x)=\overline{\eta}$ for $x\ge L$. This leads us to define the unbounded operator $B_{\theta}$ of $\mL^2(\Om)$ such that
\begin{equation}\label{defOpConjugatedPMLs}
B_{\theta}w_{\theta}=-\cfrac{1}{\gamma}\left(\beta_{\theta}\frac{\partial}{\partial x}\Big(\beta_{\theta}\frac{\partial w_{\theta}}{\partial x}\Big)+\frac{\partial^2 w_{\theta}}{\partial y^2}\right)
\end{equation}
with homogeneous Neumann boundary conditions. As $A_{\theta}$, the operator $B_{\theta}$ is not selfadjoint. Its spectral properties are summarized in the following theorem which is proved in the last section of this article.
\begin{theorem}\label{thmConjugatedPMLs} 
i) There holds 
\begin{equation}\label{essConjPMLs}
\sigma_{\mrm{ess}}(B_{\theta})=\bigcup_{n\in\N,\,t\ge0}\{n^2\pi^2+te^{-2i\theta},\,n^2\pi^2+te^{+2i\theta}\}.
\end{equation}
ii) The spectrum of $B_{\theta}$ satisfies $\sigma(B_{\theta})\subset\mathscr{R}_{\mrm{\theta}}$ with 
\begin{equation}\label{DefCPML}
\mathscr{R}_{\mrm{\theta}}:=\{z\in\Cplx\,|\,-2\theta\le\mrm{arg}(z)\le 2\theta\}.
\end{equation}
iii) Assume that $k^2\in\sigma(B_{\theta})\setminus\sigma_{\mrm{ess}}(B_{\theta})$. Then $k^2$ is real if and only if $k^2\in\mathscr{K}_{\mrm{t}}\cup\mathscr{K}_{\mrm{r}}$. Moreover if $w_{\theta}$ is an eigenfunction associated to $k^2$ such that $\pm\Im m\,k^2<0$, then $w_{\theta}\circ\mathcal{J}_{-\theta}$ is a solution of (\ref{PbInitial}) whose amplitude is exponentially growing at $\pm\infty$ and exponentially decaying at $\mp\infty$.
\end{theorem} 
\noindent The important result is that isolated real eigenvalues of $B_{\theta}$ correspond precisely to TMs and RMs. The following proposition provides a criterion to determine whether an eigenfunction associated to a real eigenvalue of $B_{\theta}$ is a TM or a RM. 
\begin{proposition}\label{PropositionCarac}
Assume that $(k^2,w_{\theta})\in\R\times \mL^2(\Om)$ is an eigenpair of $B_{\theta}$ such that $k\in(N\pi;(N+1)\pi)$, $N\in\N$. Set
\begin{equation}\label{DefAlephFunction}
\Indicator(w_{\theta})=\sum_{n=0}^N\Big|\int_{0}^1w_{\theta}(-L,y)\varphi_n(y)\,dy\Big|^2
\end{equation}
where $\varphi_n$ is defined in (\ref{defModes}). If $\Indicator(w_{\theta})=0$ then $w_{\theta}\circ\mathcal{J}_{-\theta}$ is a TM ($k^2\in\mathscr{K}_{\mrm{t}}$). If $\Indicator(w_{\theta})>0$ then $w_{\theta}\circ\mathcal{J}_{-\theta}$ is a RM ($k^2\in\mathscr{K}_{\mrm{r}}$). In this case, the incident field $u_i$ defined in (\ref{DefIncidentField}) with
\[
a_n=\int_{0}^1w_{\theta}(-L,y)\varphi_n(y)\,dy,\quad n=0,\dots,N,
\]
yields a scattered field which decays exponentially for $x\le -L$.
\end{proposition}
\noindent The next proposition tells  that  $B_{\theta}$ satisfies the celebrated $\mathcal{PT}$ symmetry property when the obstacle is symmetric with respect to the $(Oy)$ axis. This ensures in particular the stability of simple real eigenvalues, with respect to perturbations of the obstacle satisfying the same symmetry constraint. 
\begin{proposition}\label{PropositionPTsym}
Assume that $\gamma$ satisfies $\gamma(x,y)=\gamma(-x,y)$ for all $(x,y)\in\Om$. Then the operator $B_{\theta}$ is $\mathcal{PT}$-symmetric ($\mathcal{PT}B_{\theta}\mathcal{PT}=B_{\theta}$) with $\mathcal{P}\varphi(x,y)=\varphi(-x,y)$, $\mathcal{T}\varphi(x,y)=\overline{\varphi(x,y)}$ for $\varphi\in\mL^2(\Om)$. Therefore, we have $\sigma(B_{\theta})=\overline{\sigma(B_{\theta})}$.
\end{proposition}
\noindent The proof is straightforward observing that the $\beta_{\theta}$ defined after (\ref{DefPbFortCPMLs}) satisfies $\beta_{\theta}(-x,y)=\overline{\beta_{\theta}(x,y)}$.\\
\newline
Finally let us mention a specific difficulty which appears in the spectral analysis of $B_{\theta}$. While Theorem \ref{thmUsualPMLs} guarantees that $\sigma(A_{\theta})\setminus\sigma_{\mrm{ess}}(A_{\theta})$ is discrete, we do not write such a statement for the operator $B_{\theta}$ in Theorem \ref{thmConjugatedPMLs}. A major difference between both operators is that $\Cplx\setminus \sigma_{\mrm{ess}}(A_{\theta})$ is connected whereas $\Cplx\setminus \sigma_{\mrm{ess}}(B_{\theta})$ has a countably infinite number of connected components. As a consequence, to prove that $\sigma (B_{\theta})\setminus \sigma_{\mrm{ess}}(B_{\theta})$ is discrete using the Fredholm analytic theorem, it is necessary to find one $\lambda$ such that $B_{\theta}-\lambda$ is invertible in each of the components of $\Cplx\setminus \sigma_{\mrm{ess}}(B_{\theta})$. In general, in presence of an obstacle, such a $\lambda$ probably exists (proofs for certain classes of $\gamma$ can be obtained working as in \cite{BoCN15}). But for this problem, we can have surprising perturbation results. Thus, if there is no obstacle ($\gamma\equiv1$ in $\Om$), then there holds $\sigma(B_{\theta})=\mathscr{R}_{\mrm{\theta}}$ (see (\ref{DefCPML})): all connected components of $\Cplx\setminus \sigma_{\mrm{ess}}(B_{\theta})$, except the one containing the  complex half-plane $\Re e\,\lambda<0$, are filled with eigenvalues. 
 To show this result, observe that for $k^2\in\sigma(B_{\theta})\setminus\sigma_{\mrm{ess}}(B_{\theta})$, the function $u\circ \mathcal{J}_{\theta}$, with $u(x,y)=e^{ik x}$, is a non-zero element of $\ker B_{\theta}$. Notice that this pathological property is also true when $\Om$ contains a family of sound hard cracks (homogeneous Neumann boundary condition) parallel to the $(Ox)$ axis.

\section{Numerical experiments}\label{SectionNumerics}

\subsection{ Classical complex scaling: classical complex resonance modes} 
\begin{figure}[!ht]
\centering
\includegraphics[width=.49\linewidth]{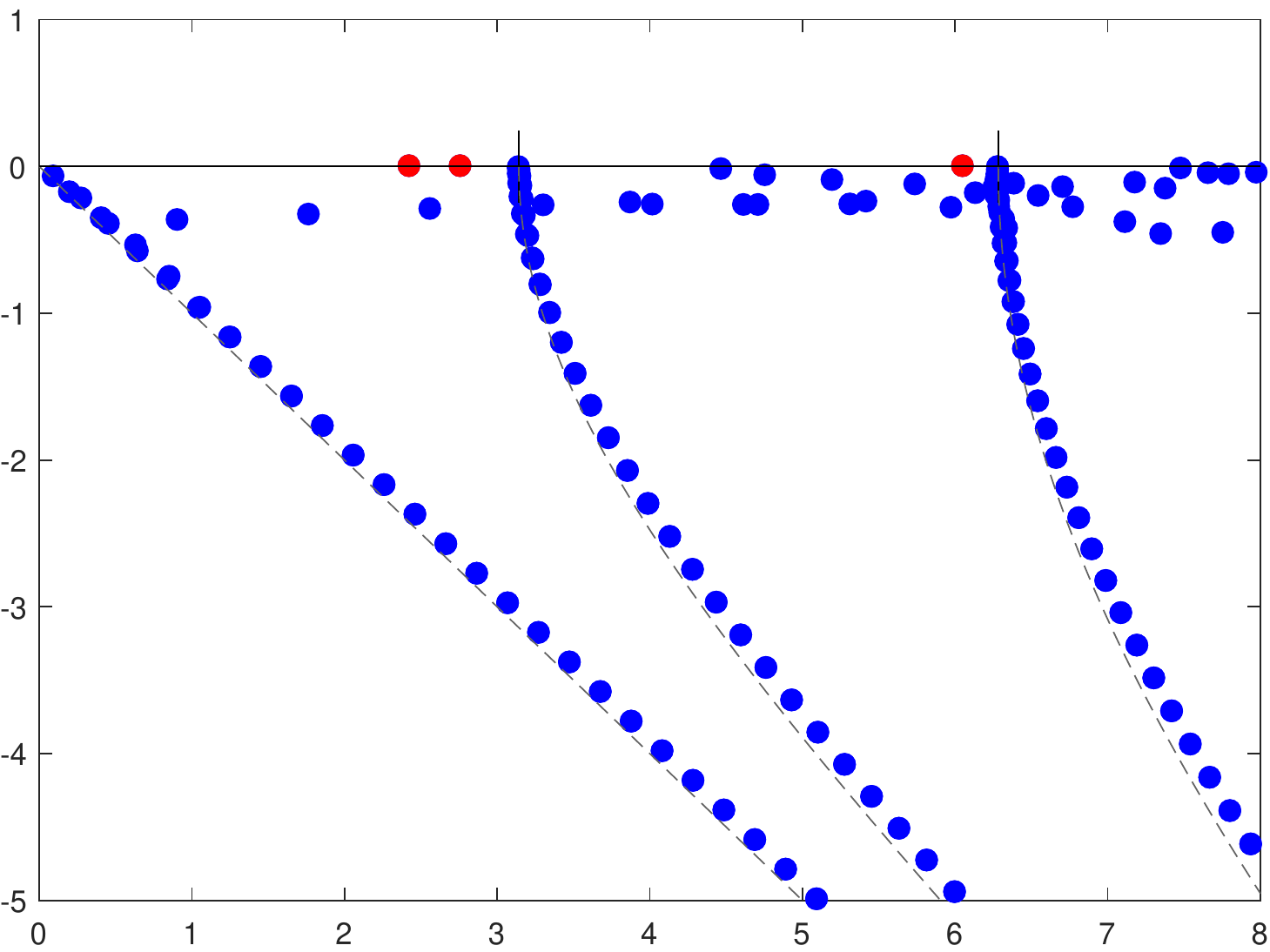}\ \includegraphics[width=.49\linewidth]{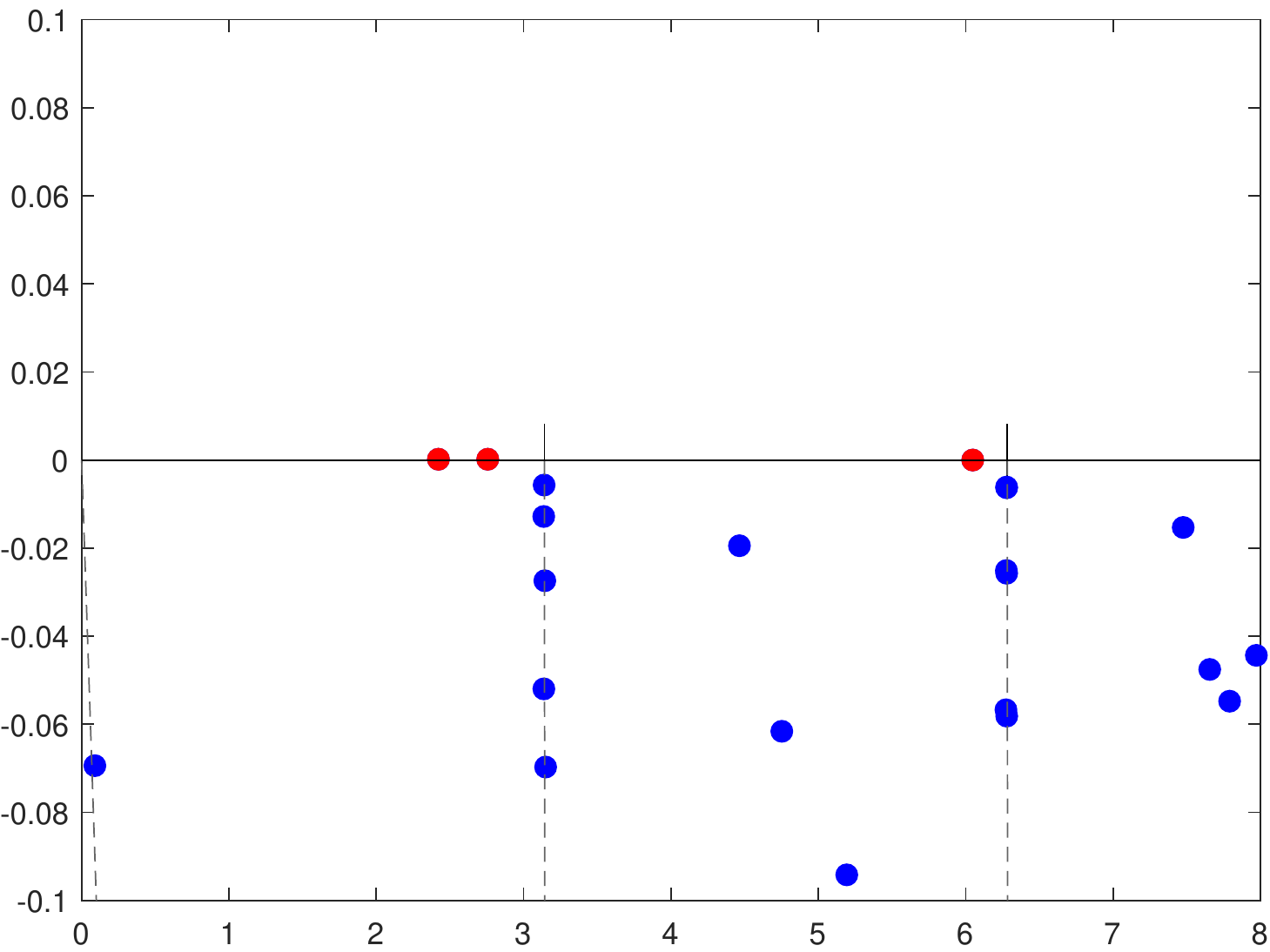}
\caption{Classical complex resonances in the complex $k$ plane corresponding to the spectrum of $A_{\theta}$ for a symmetric obstacle (Figure \ref{figsetting} (a)). The trapped modes are in red, the dashed lines represent the essential spectrum of  $A_{\theta}$ (see (\ref{essClassicalPMLs})). The picture on the right is a zoom-in of that on the left.} 
\label{figClassicalPMLs}
\end{figure}

We first compute the spectrum of the operator $A_{\theta}$ defined in (\ref{defOpClassicalPMLs}) with a classical complex scaling (complex resonance spectrum).  For the numerical experiments, we truncate the computational domain at some distance of the obstacle and use finite elements. This corresponds to the so-called \textit{Perfectly Matched Layers} (PMLs) method. We refer the reader to \cite{Kalv13} for the numerical analysis of the error due to truncation of the waveguide and discretization. The setting is as follows. We take $\gamma$ such that $\gamma=5$ in $\mathscr{O}=(-1;1)\times(0.25;0.75)$ and $\gamma=1$ in $\Om\setminus\overline{\mathscr{O}}$ (see Figure \ref{figsetting} (a)).  In the definition of the maps $\mathcal{I}_{\theta}$, $\alpha_{\theta}$ (see (\ref{defDilatationx}), (\ref{defEquationDilatation})), we take $\theta=\pi/4$ (so that $\eta=e^{i\pi/4}$) and $L=1$. In practice, we use a $\mrm{P2}$ finite element method in the bounded domain $\Om_{12}=\{(x,y)\in\Om\,|\,-12<x<12\}$ with Dirichlet boundary condition at $x=\pm12$. Finite element matrices are constructed with \textit{FreeFem++}\footnote{\textit{FreeFem++}, \url{http://www.freefem.org/ff++/}.}. \\
\newline
In Figure \ref{figClassicalPMLs} and in the rest of the paper, we display the square root of the spectrum ($k$ instead of $k^2$). The vertical marks on the real axis correspond to the thresholds ($0$, $\pi$, $2\pi$, ...). In accordance with Theorem \ref{thmUsualPMLs}, we observe that $\sqrt{\sigma(A_{\theta})}$ is located in the region $\sqrt{\mathscr{R}_{\theta}^-}=\{z\in\Cplx\,|\,-\theta\le\mrm{arg}(z)\le 0\}$. Moreover, the discretisation of the essential spectrum $\sigma_{\mrm{ess}}(A_{\theta})$ defined in (\ref{essClassicalPMLs}) and forming branches starting at the threshold points appears clearly. Note that a simple calculation shows that $\sqrt{\{n^2\pi^2+te^{-2i\theta},\,t\ge0\}}$ is a half-line for $n=0$ and a piece of hyperbola for $n\ge1$. This is precisely what we get. 
Eigenvalues located on the real axis correspond to trapped modes ($k^2\in\mathscr{K}_{\mrm{t}}$).
 In the chosen setting, which is symmetric with respect to the axis $\R\times\{0.5\}$, one can prove that trapped modes exist \cite{EvLV94}. 
On the other hand, the eigenvalues in the complex plane which are not the discretisation of the essential spectrum correspond to complex resonances. 

\subsection{Conjugated complex scaling: reflectionless modes} 
\begin{figure}
\centering
\includegraphics[width=.49\linewidth]{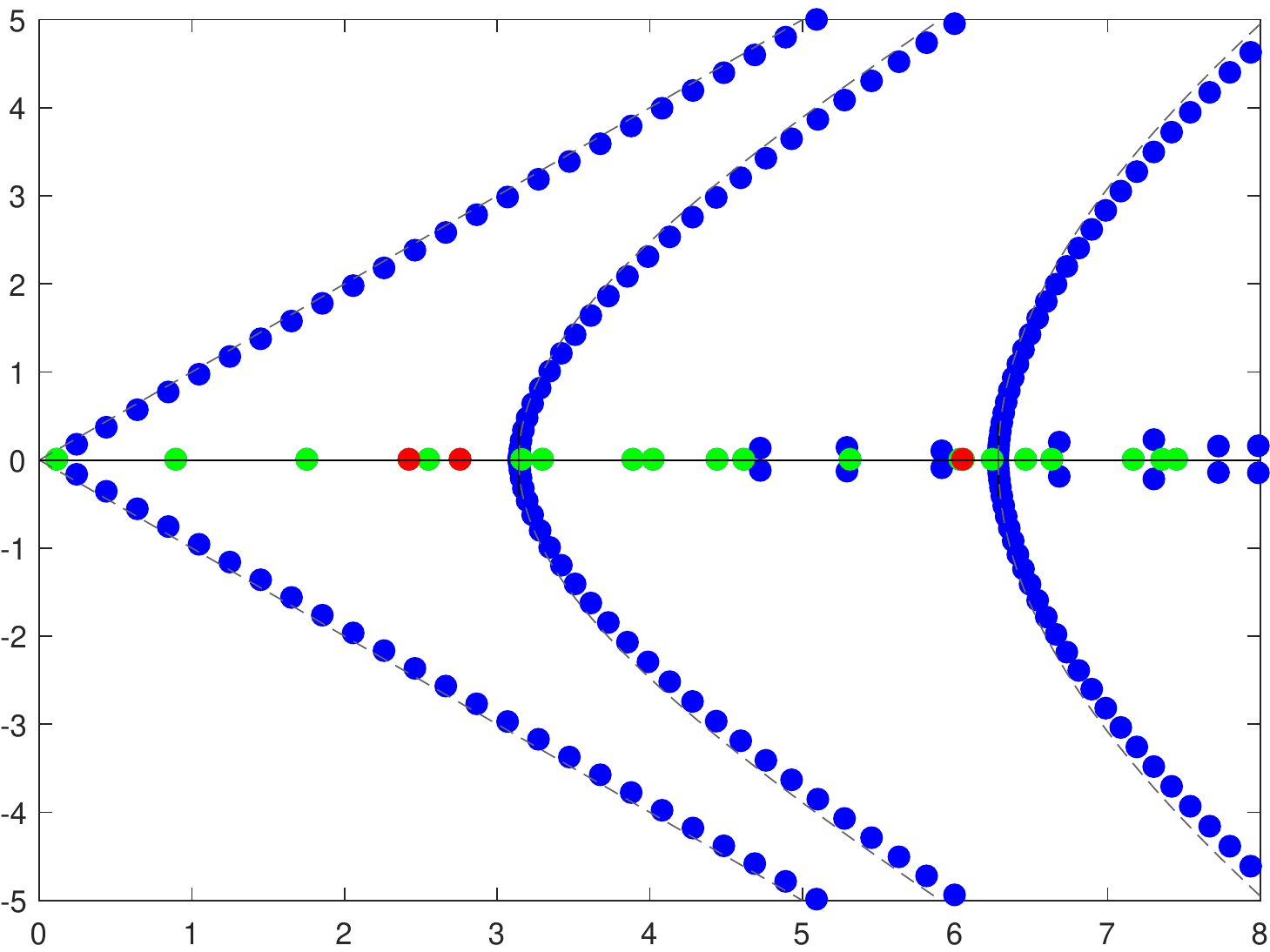}\ \includegraphics[width=.4805\linewidth]{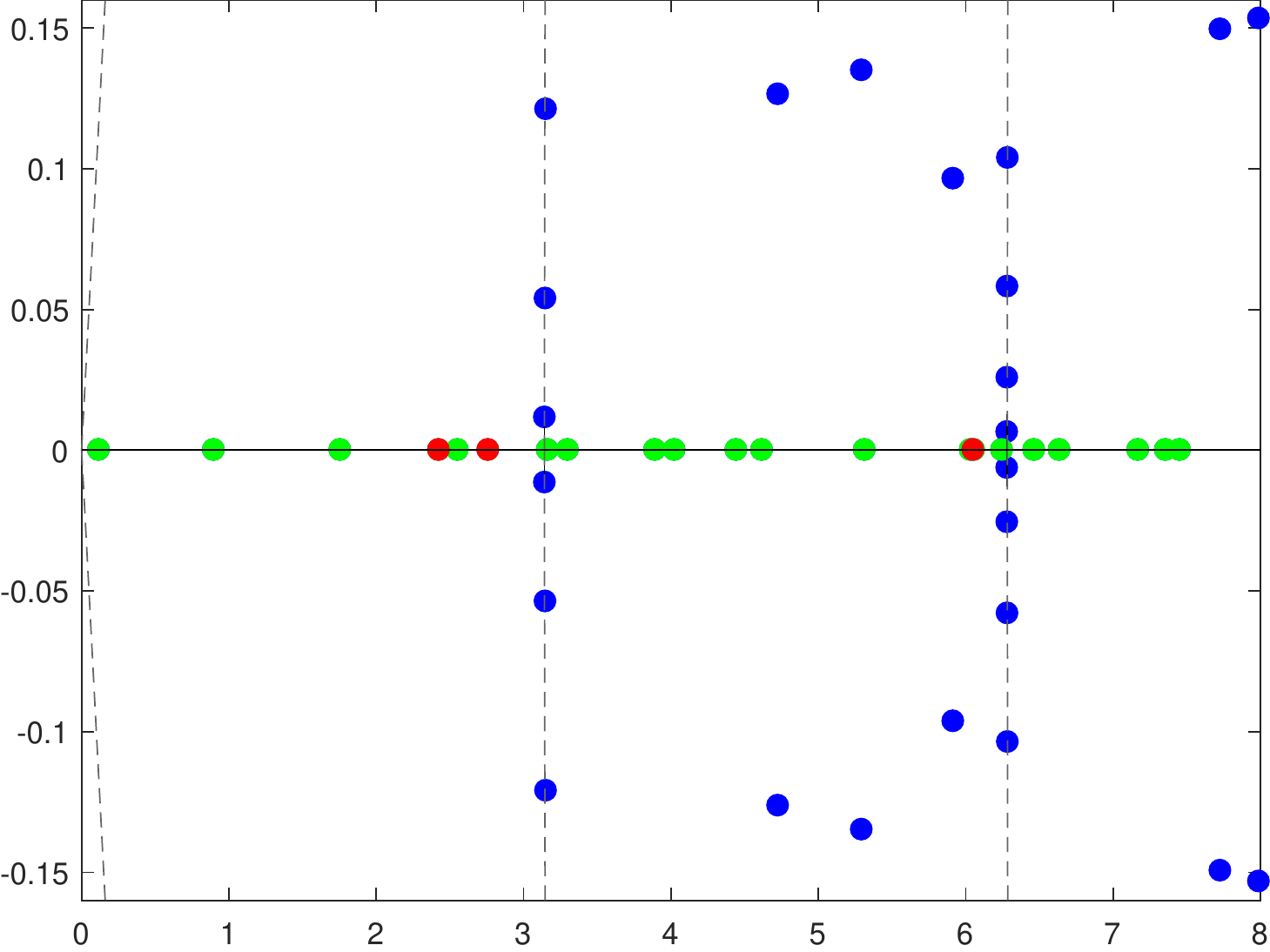}
\caption{Reflectionless eigenvalues in the complex $k$ plane corresponding to the spectrum of $B_{\theta}$ for a symmetric obstacle (Figure \ref{figsetting} (a)). The trapped modes in red are the same as in Figure \ref{figClassicalPMLs}. The reflectionless modes are in green and the dashed lines represent the essential spectrum of  $B_{\theta}$ (see (\ref{essConjPMLs})). The picture on the right is a zoom-in of that on the left.
}
\label{figConjugatedPMLs}
\end{figure}

\begin{figure}
\centering
\includegraphics[width=.49\linewidth]{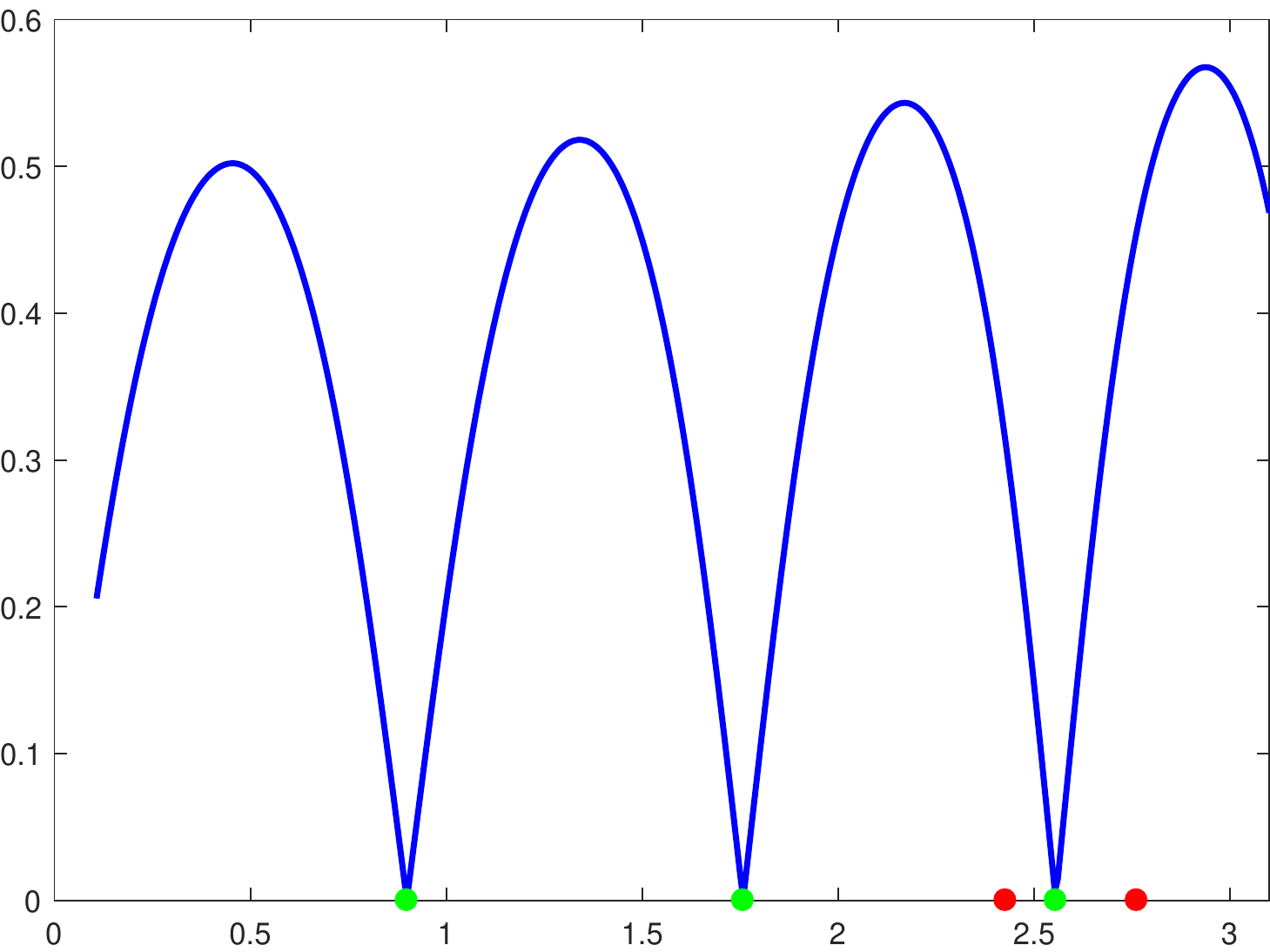}
\caption{Curve $k\mapsto |R_{00}(k)|$ (modulus of the reflection coefficient) for $k\in(0.1;3.1)$. The green and red dots represent respectively the reflectionless modes and the trapped modes computed in Figure \ref{figConjugatedPMLs}. We indeed observe that $R_{00}(k)$ is null for reflectionless $k$.}
\label{figBalayage}
\end{figure}

\begin{figure}[!ht]
\centering
\includegraphics[width=.65\linewidth]{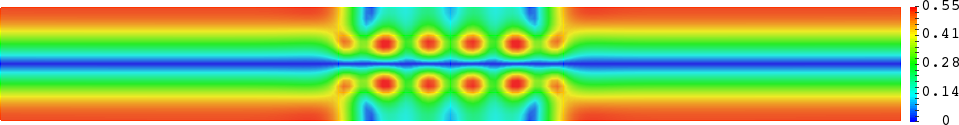}\\[1pt]
\includegraphics[width=.65\linewidth]{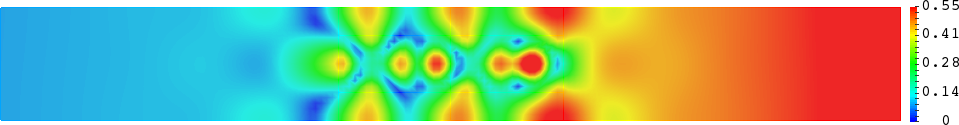}\\[1pt]
\includegraphics[width=.65\linewidth]{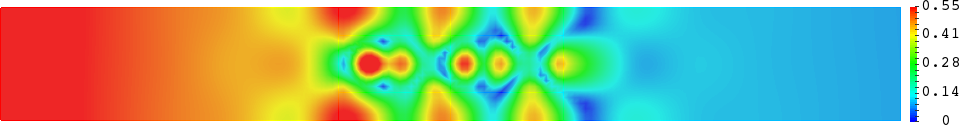}
\caption{Modulus of eigenfunctions of $B_{\theta}$ associated to the eigenvalues $k\approx 5.31$ (top), $k\approx 5.29 - 0.13i$ (middle) and $k\approx 5.29+0.13i$ (bottom) obtained in Figure \ref{figConjugatedPMLs}. The symmetry $x\to-x$ for real $k$ (due to $\mathcal{PT}$-symmetry) disappears for complex $k$. }
\label{figBrokenSym}
\end{figure}

Now we compute the spectrum of the operator $B_{\theta}$ defined in (\ref{defOpConjugatedPMLs}) with a conjugated complex scaling. First, we use exactly the same symmetric setting (see Figure \ref{figsetting} (a)) as in the previous paragraph. In Figure \ref{figConjugatedPMLs}, we display the square root of the spectrum $\sqrt{\sigma(B_{\theta})}$. Since $\gamma$ satisfies $\gamma(x,y)=\gamma(-x,y)$, according to Proposition \ref{PropositionPTsym} we know that $B_{\theta}$ is $\mathcal{PT}$-symmetric and that therefore its spectrum is stable by conjugation ($\sigma(B_{\theta})=\overline{\sigma(B_{\theta})}$). This is indeed what we obtain. Note that the mesh has been constructed so that $\mathcal{PT}$-symmetry is preserved at the discrete level. $\mathcal{PT}$-symmetry is an interesting property in our case because it guarantees that eigenvalues located close to the real axis which are isolated (no other eigenvalue in a vicinity) are real. Therefore, according to Theorem \ref{thmConjugatedPMLs}, they correspond to trapped modes or to reflectionless modes. 
Remark that, for the same geometry, the spectrum of $B_{\theta}$ (Figure \ref{figConjugatedPMLs}) contains more elements on the real axis than the spectrum of $A_{\theta}$ (Figure \ref{figClassicalPMLs}): the additional elements (green points in Figure \ref{figConjugatedPMLs}) correspond to reflectionless modes. 
\begin{figure}
\centering
\includegraphics[width=.65\linewidth]{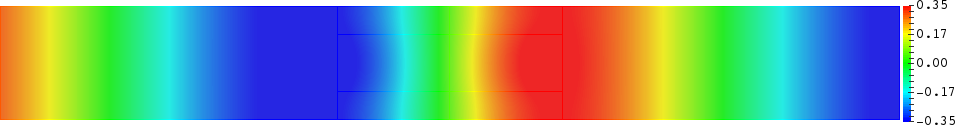}\\[1pt]
\includegraphics[width=.65\linewidth]{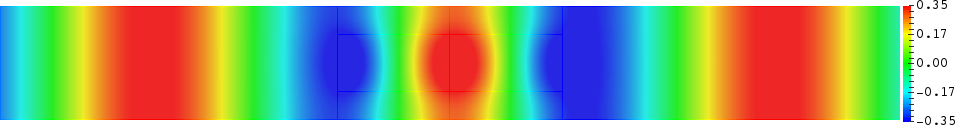}\\[1pt]
\includegraphics[width=.65\linewidth]{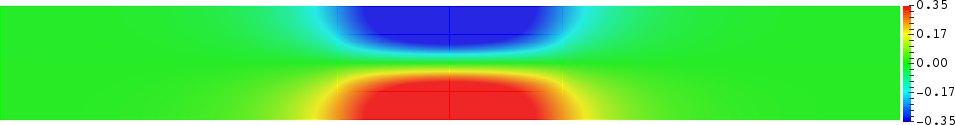}\\[1pt]
\includegraphics[width=.65\linewidth]{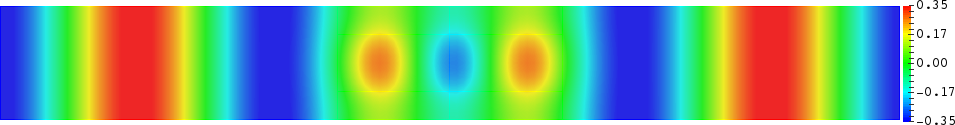}\\[1pt]\includegraphics[width=.65\linewidth]{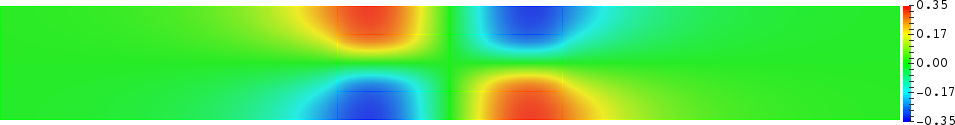}\\[1pt]
\includegraphics[width=.65\linewidth]{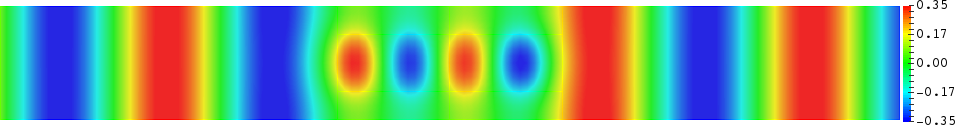}\\[1pt]
\includegraphics[width=.65\linewidth]{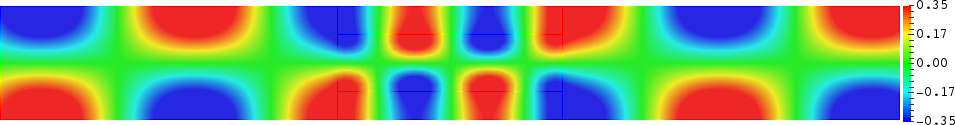}\\[10pt]
\renewcommand{\arraystretch}{1.8}
\begin{tabular}{|>{\centering\arraybackslash}m{1.2cm}|>{\centering\arraybackslash}m{0.7cm}|>{\centering\arraybackslash}m{0.7cm}|>{\centering\arraybackslash}m{1.5cm}|>{\centering\arraybackslash}m{0.7cm}|>{\centering\arraybackslash}m{1.5cm}|>{\centering\arraybackslash}m{0.7cm}|>{\centering\arraybackslash}m{0.7cm}|}
\hline 
\cellcolor{Gray} $k$ &   0.9 &  1.8 &  2.4 &  2.6 &  2.8 &  3.3 &  3.9 \\
\hline
\cellcolor{Gray} $\Indicator(w_{\theta})$  & 0.14 & 0.14 & 8.0\,$10^{-10}$ & 0.14 & 4.3\,$10^{-9}$ & 0.14& 0.14 \\
\hline
\end{tabular}
\caption{Top: real part of eigenmodes associated with real eigenvalues of $B_{\theta}$ from Figure \ref{figConjugatedPMLs}. 
Bottom: value of $k$ and of the indicator function $\Indicator$ for each of these $7$ eigenmodes. The 3rd and the 5th eigenmodes are trapped modes, the five others are reflectionless modes.}
\label{figEigenfunctions}
\end{figure}
In Figure \ref{figEigenfunctions} top, we represent the real part of eigenfunctions associated with seven real eigenvalues of $B_{\theta}$. 
To obtain these pictures, we take $L=4$ in the definition of $\mathcal{J}_{\theta}$ in (\ref{defDilatationConjugated})
and we display only the restrictions of the eigenfunctions to $\Om_{L}=\{(x,y)\in\Om\,|\,-L<x<L\}$. 
We recognize two trapped modes (images 3 and 5). 
The other modes are reflectionless modes. In Figure \ref{figEigenfunctions} bottom, we provide the value of the indicator function $\Indicator$ defined in (\ref{DefAlephFunction}) for the seven eigenmodes. We have to mention that eigenmodes are normalized so that their $\mL^2$ norm  is equal to one. The indicator function $\Indicator$ offers a clear criterion to distinguish between trapped modes and reflectionless modes. 
Moreover, in order to inspect the scattering coefficient, we remark that for reflectionless modes associated with wavenumbers $k$ smaller than $\pi$, the incident field $u_i$ in (\ref{DefIncidentField}) decomposes only on the piston mode $w^{+}_0(x,y)=e^{ikx}/\sqrt{2k}$ (monomode regime). In this case the reflection matrix $R(k)$ in (\ref{reflection matrix}) is nothing but the usual reflection coefficient. 
In Figure \ref{figBalayage}, we thus display the modulus of this coefficient $R_{00}(k)$ with respect to $k\in(0.1;3.1)$. As expected, we observe that $R_{00}$ vanishes for the values of $k$ obtained in Figure \ref{figEigenfunctions} solving the spectral problem for $B_{\theta}$. 
Of course obtaining the curve $k\mapsto |R_{00}(k)|$ is relatively costly and it is precisely 
what we want to avoid by computing the reflectionless $k$ as eigenvalues. Here it is simply a way to check our results. 
In Figure \ref{figBrokenSym}, we represent the modulus of reflectionless mode eigenfunctions of $B_{\theta}$ associated with one real eigenvalue and two complex conjugated eigenvalues. We observe, and this is true in general, a symmetry with respect to the $(Oy)$ axis for modes corresponding to real eigenvalues which disappears for complex ones. This is the so-called broken symmetry phenomenon which is well-known for $\mathcal{PT}$-symmetric operators (see e.g. the review \cite{Bend07}).

\begin{figure}[!ht]
\centering
\includegraphics[width=.49\linewidth]{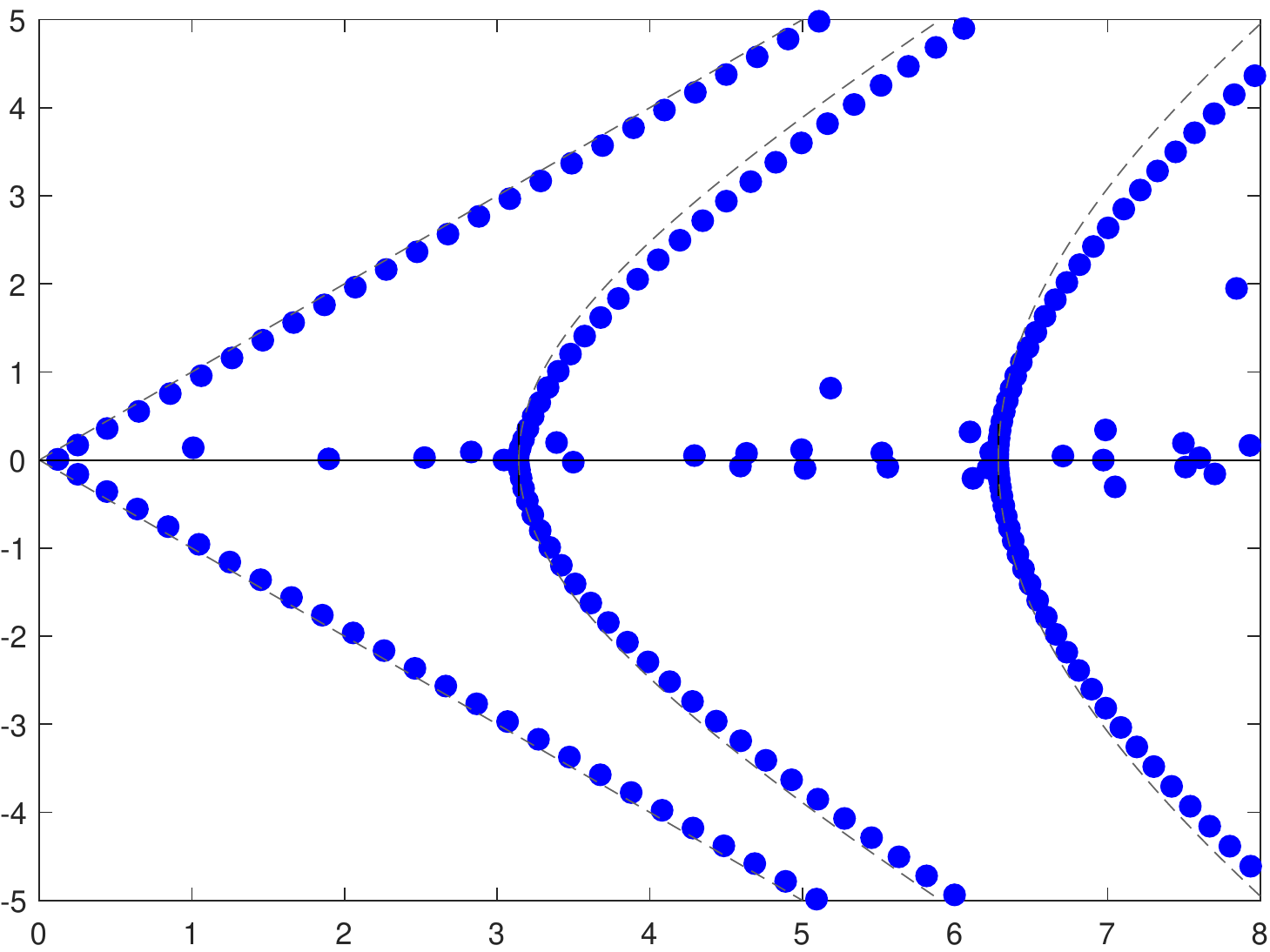}\ \includegraphics[width=.4805\linewidth]{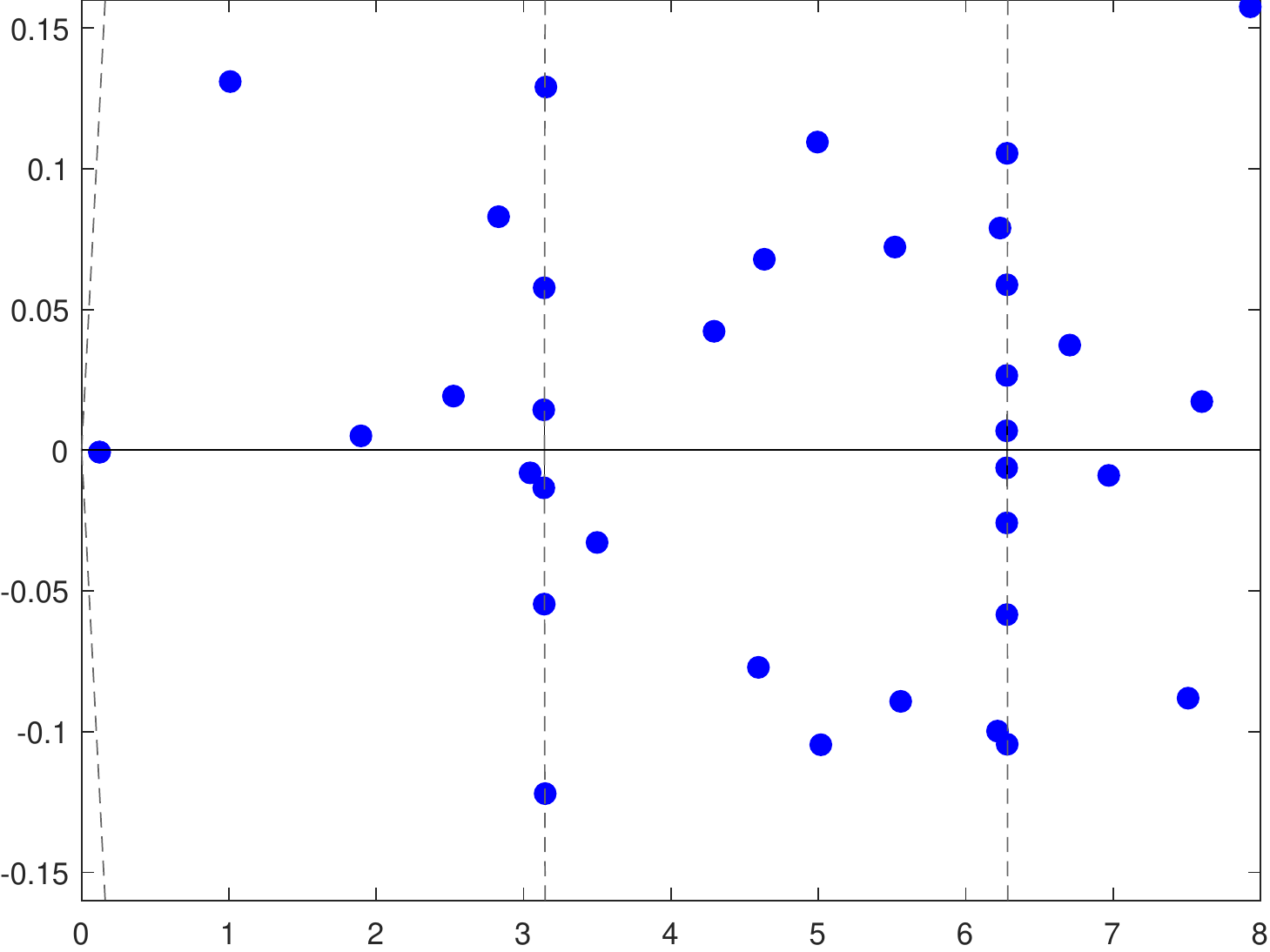}
\caption{Spectrum of $B_{\theta}$ in the complex $k$ plane for a non symmetric obstacle (Figure \ref{figsetting} (b)). The dashed lines represent the essential spectrum of  $B_{\theta}$ (see (\ref{essConjPMLs})). The spectrum is not stable by conjugation. The picture on the right is a zoom-in of that on the left.}
\label{figSpectrumNonSym}
\end{figure}

\begin{figure}[!ht]
\centering 
\includegraphics[width=.49\linewidth]{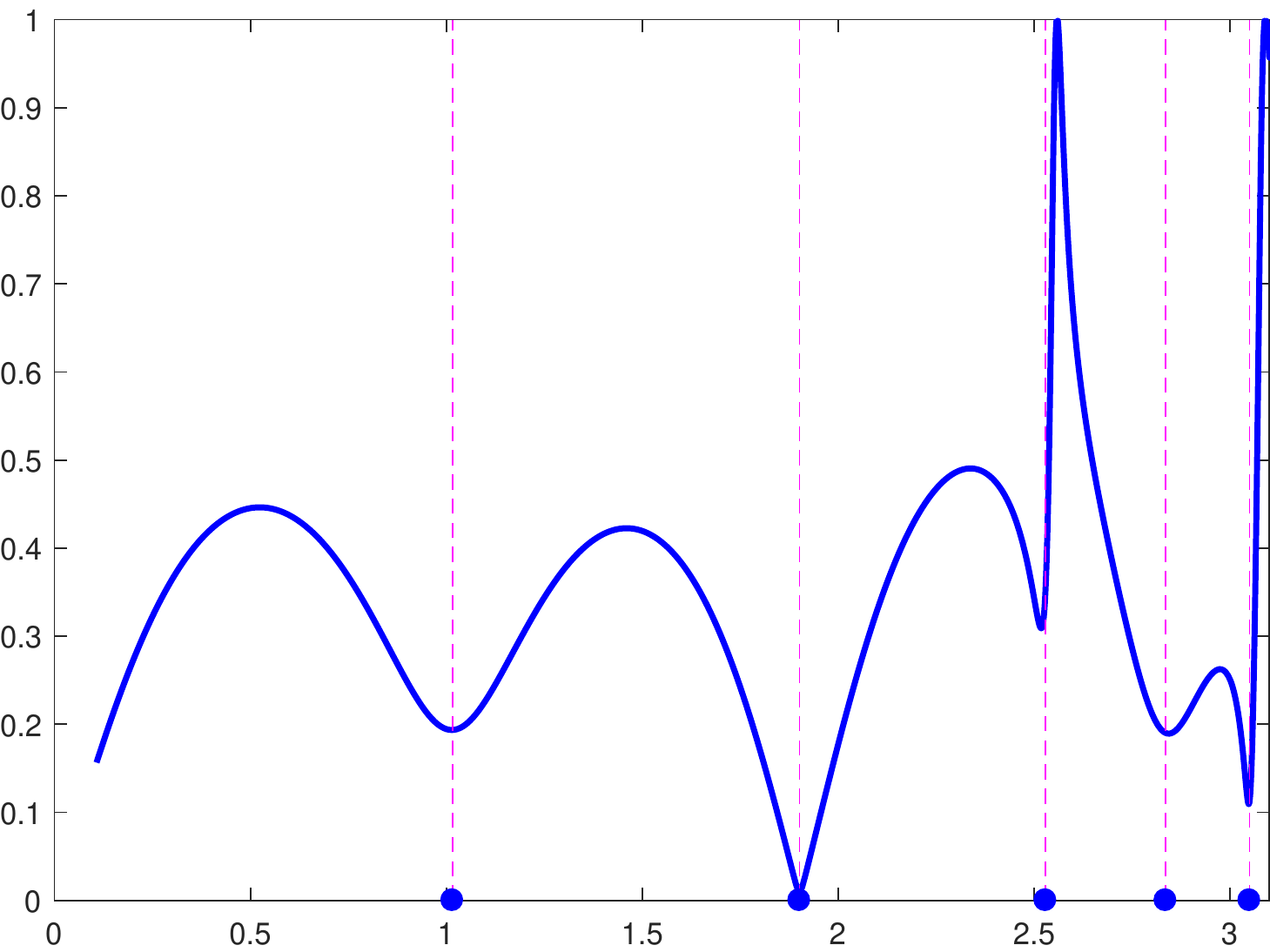}
\caption{Curve $k\mapsto |R_{00}(k)|$ (modulus of the reflection coefficient) for $k\in(0.1;3.1)$ and a non symmetric obstacle. The blue dots and the vertical dashed lines correspond to the real parts of the eigenvalues of $B_{\theta}$ located close to the real axis computed in Figure \ref{figSpectrumNonSym}. We observe that $|R_{00}(k)|$ is minimal for these particular $k$.}
\label{figNonSymCurve}
\end{figure}

\noindent Now, we use the  non-symmetric setting (see Figure \ref{figsetting} (b)), 
and we display the square root of the spectrum of $B_{\theta}$ (in Figure \ref{figSpectrumNonSym})
for a coefficient $\gamma$ which is not symmetric in $x$ nor in $y$. 
More precisely, we take $\gamma$ such that $\gamma=5$ in $\mathscr{O}=(-1;0]\times(0.25;0.5)\cup[0;1)\times(0.25;0.75)$ and $\gamma=1$ in $\Om\setminus\overline{\mathscr{O}}$. We observe that the  spectrum is no longer stable by conjugation ($\sigma (B_{\theta})\ne\overline{\sigma (B_{\theta})}$) since
the operator $B_{\theta}$ is not  $\mathcal{PT}$-symmetric, and there is no ``help'' for the eigenvalues to be real.
However, a closer look shows the presence of eigenvalues close to the real axis,
in particular for $k\approx1.0 + 0.13i$, $k\approx1.9 + 0.005i$, $k\approx2.5 + 0.02i$, $k\approx2.8 + 0.08i$ and $k\approx3.0 - 0.008i$. In Figure \ref{figNonSymCurve}, we represent $k\mapsto|R_{00}(k)|$ for $k\in(0.1;3.1)$ where there is only one propagating mode in the leads.
 It is interesting to note that the above computed complex reflectionless modes (located close to the real axis)   have an influence on this curve. More precisely, $k\mapsto|R_{00}(k)|$ attains minima for $k\in(0.1;3.1)$ close to the real part of these complex reflectionless modes. Therefore complex reflectionless modes also have  significance for scattering at real frequncies.

\newpage

\section{Proofs}\label{SectionProofs}
Finally we give the proofs of Theorem \ref{thmConjugatedPMLs} and Proposition \ref{PropositionCarac}.\\
\newline
\textbf{Proof of Theorem \ref{thmConjugatedPMLs}.} 
$i)$ First, let us explain how to show that for $t\ge0$ and $n\in\N$, $n^2\pi^2+te^{-2i\theta}$ belongs to $\sigma_{\mrm{ess}}(B_{\theta})$. 
Consider a smooth cut-off function $\chi$ defined on $\R$ such that $\chi(x)=0$ for $|x|>L$ and $\|\chi\|_{\mL^2(\R)}=1$. For $m\ge1$, set $v_+^{(m)}(x,y)=m^{-1/2}\chi\left((x-m^2L)/m\right)e^{i\sqrt{t}x}\varphi_n(y)$ where $\varphi_n$ is defined in (\ref{defModes}). One can check that $(v_+^{(m)})$ is a singular sequence for $B_{\theta}$ at $n^2\pi^2+te^{-2i\theta}$ (note that the supports of $v_+^{(m)}$ and $v_+^{(m')}$ for $m\neq m'$ do not overlap). On the other hand, we prove that $(v_-^{(m)})$, with $v_-^{(m)}(x,y)=v_+^{(m)}(-x,y)$, is a singular sequence
for $B_{\theta}$ at $n^2\pi^2+te^{+2i\theta}$. Summing up, we obtain
$\cup_{n\in\N,\,t\ge0}\{n^2\pi^2+te^{-2i\theta},\,n^2\pi^2+te^{+2i\theta}\}\subset \sigma_{\mrm{ess}}(B_{\theta})$.\\
The converse inclusion requires a bit more work. Observe that if 
$\lambda\in\sigma_{\mrm{ess}}(B_{\theta})$, then one can prove by localization the existence of a corresponding singular sequence $(v^{(m)})$ supported either in $x>L$ or in $x<-L$. In the first case, it means that $(v^{(m)})$ is also a singular sequence for $A_{\theta}$, so that $\lambda\in\sigma_{\mrm{ess}}(A_{\theta})$. In the second case, $(v^{(m)})$  is a singular sequence for the complex conjugate of $A_{\theta}$, which implies that $\overline{\lambda}\in\sigma_{\mrm{ess}}(A_{\theta})$. Finally, the result follows from item $i)$ of Theorem \ref{thmUsualPMLs}. 
\\[4pt]
$ii)$ Define the form $b_{\theta}(\cdot,\cdot)$   such that 
\[
b_{\theta}(u,v)=\int_{\Om}\beta_{\theta}\partial_xu\,\partial_x\overline{v}+\beta^{-1}_{\theta}\partial_yu\,\partial_y\overline{v}-k^2\,\beta^{-1}_{\theta}u\,\overline{v}\,dxdy.
\]
One can check that for $k^2\in\Cplx\setminus\mathscr{R}_{\theta}$, the numbers $1$, $e^{\pm i\theta}$, $-k^2$, $-k^2\,e^{\pm i\theta}$ are all located in some region $\{z\in\Cplx\,|\,\Re e\,(e^{i\kappa}z)\ge \mu\}$ for some $\kappa\in [0;2\pi)$ and $\mu>0$. We deduce that 
\[
\Re e\,(e^{i\kappa}b_{\theta}(u,u))\ge \mu \int_{\Om}|u|^2+|\partial_xu|^2+|\partial_yu|^2\,dxdy
\]
which shows that $b_{\theta}(\cdot,\cdot)$ is coercive on the Sobolev space $\mH^1(\Omega)$. As a consequence of the Lax-Milgram theorem, $B_{\theta}-k^2$ is invertible for $k^2\in\Cplx\setminus\mathscr{R}_{\theta}$, which guarantees that $\sigma(B_{\theta})\subset\mathscr{R}_{\theta}$.\\[4pt]
$iii)$ If $k^2\in\mathscr{K}_{\mrm{t}}\cup\mathscr{K}_{\mrm{r}}$, then  $k^2$ is real and by construction $k^2\in\sigma(B_{\theta})\setminus\sigma_{\mrm{ess}}(B_{\theta})$. Conversely, assume that $k^2$ is real and that $k^2\in\sigma(B_{\theta})\setminus\sigma_{\mrm{ess}}(B_{\theta})$. There is one $N\in\N$ such that $k\in(N\pi;(N+1)\pi)$. Consider a non zero  $w_{\theta}\in\ker (B_{\theta}-k^2)$. Setting $u=w_{\theta}\circ\mathcal{J}_{-\theta}$, we find that $u$ satisfies $\Delta u+k^2\gamma u=0$ in $\Om$ and expands as
\begin{equation}\label{DecompoEigenBtheta}
u=\begin{array}{|l}
\dsp\sum_{n=0}^{N} a^{-}_n\,w_n^{+}+\sum_{n=N+1}^{+\infty} a^{-}_n\,w_n^{-}\quad\mbox{ for } x\le -L\\[12pt]
\dsp\sum_{n=0}^{+\infty} a^{+}_n\,w_n^{+}\quad\mbox{ for } x\ge L,
\end{array}
\end{equation}
with $(a^{\pm}_n)\in\Cplx^{\N}$. If one of the $a^{-}_n$, for $n=0,\dots,N$, is non zero, it means that $u$ is a RM. If on the contrary $a^{-}_n=0$ for all $n=0,\dots,N$, then one can prove that $a^{+}_n=0$ for all $n=0,\dots,N$, so that $u$ is a TM. Indeed, integrating by parts, we show that the quantity
\[
\mathscr{F}=\int_{\Sigma_{-L}\cup\Sigma_L} (\partial_x u \overline{u}-u \partial_x \overline{u}) \,dy,
\]
with $\Sigma_{\pm L}=\{\pm L\}\times(0;1)$, satisfies $\mathscr{F}=0$. On the other hand, a direct calculation using expansion (\ref{DecompoEigenBtheta}) and the orthonormality of the $\varphi_n$ yields 
$
\mathscr{F}=\sum_{n=0}^N i(|a_n^+|^2+|a_n^-|^2)
$
and the result follows.
\\
Finally, consider some $k^2\in\sigma(B_{\theta})\setminus\sigma_{\mrm{ess}}(B_{\theta})$ such that $\Im m\,k^2>0$ (the case $\Im m\,k^2<0$ is similar). There is a unique $N\in\N$ such that $2\pi-2\theta<\mrm{arg}(k^2-(N\pi)^2)<2\pi$ and $\pi<\mrm{arg}(k^2-((N+1)\pi)^2)<2\pi-2\theta$. Then if $w_{\theta}\in\ker (B_{\theta}-k^2)\setminus\{0\} $, expansion (\ref{DecompoEigenBtheta}) for $u=w_{\theta}\circ\mathcal{J}_{-\theta}$ holds and one of the $a^{-}_n$, $n=0,\dots,N$, has to be non zero. Indeed, otherwise $u$ would be exponentially decaying for $\pm x\ge L$ and $k^2$ would be in $\sigma(A)$, which is impossible because $\sigma(A)=[0;+\infty)$. Therefore the amplitude of $u$ is exponentially decaying at $+\infty$ and exponentially growing at $-\infty$.\hfill$\square$\\
\newline
\textbf{Proof of Proposition \ref{PropositionCarac}.} 
If $(k^2,w_{\theta})\in\R\times \mL^2(\Omega)$ is an eigenpair of $B_{\theta}$, then  $u=w_{\theta}\circ\mathcal{J}_{-\theta}$ expands as in (\ref{DecompoEigenBtheta}). Moreover, we deduce from the orthogonality of the $\varphi_n$ that 
\[
a_n^-=\int_{0}^1w_{\theta}(-L,y)\varphi_n(y)\,dy,\quad n=0,\dots,N,
\]
which gives the result, using the same arguments as in the proof of Theorem \ref{thmConjugatedPMLs}, item $iii)$. 
\hfill$\square$

\section{Concluding remarks}\label{SectionConclusion}
It is often desirable to determine frequencies for which a wave can be completely transmitted through a structure, a task usually leading to the tedious work of evaluating the scattering coefficients for each frequency. Here, we have shown that reflectionless frequencies can be directly computed as the eigenvalues of a non-selfadjoint operator $B_{\theta}$ (see (\ref{defOpConjugatedPMLs})) with conjugated complex scalings enforcing ingoing behaviour in the incident lead and outgoing behaviour in the other lead. The reflectionless spectrum of this operator $B_{\theta}$ provides a  complementary information to the one contained in the classical complex resonance spectrum associated with leaky modes which decompose only on outgoing waves (see the operator $A_\theta$ in (\ref{defOpClassicalPMLs})). Note that eigenvalues corresponding to trapped modes belong to both the reflectionless spectrum and to the classical complex resonance spectrum because trapped modes do not excite propagating waves.\\
\newline
Moreover, since ingoing and outgoing complex scalings can be associated, respectively, with gain and loss, one observes that the non-selfadjoint operator $B_{\theta}$ leads to consider a natural $\mathcal{P}\mathcal{T}$ symmetric problem when the structure has mirror symmetry. Interestingly, a direct calculus shows that in the very simple case of a $\mrm{1D}$ transmission problem through a slab of constant index, reflectionless frequencies are all real. This gives an example of a non-selfadjoint $\mathcal{P}\mathcal{T}$ symmetric operator with only real eigenvalues.\\
\newline
In this work, we investigated scattering problems in waveguides with $N=2$ leads for which two reflectionless spectra exist: one associated with incident waves propagating from the left and another corresponding to incident waves propagating from the right. The more general case with $N$ ($N \ge 2$) leads can be considered as well. Among the total of $2^N$ different spectra with an ingoing or an outgoing complex scaling in each lead, two spectra correspond to eigenmodes which decompose on waves which are all outgoing or all ingoing. As a consequence, there are $2^N-2$ reflectionless spectra.

\bibliography{Bibli}
\bibliographystyle{plain}
\end{document}